\documentclass {article}

\usepackage {amsmath}
\usepackage {amsthm}
\usepackage {amssymb}
\usepackage {enumerate}
\usepackage {xypic}

\newtheorem* {thm*} {Theorem}
\newtheorem {thm}   {Theorem}
\newtheorem {prop}  {Proposition}
\newtheorem {lem}   {Lemma}
\newtheorem {cor}   {Corollary}

\newtheorem* {conj*}{Conjecture}
\newtheorem* {convention} {Convention}
\newtheorem* {claim}{Claim}

\newtheoremstyle {numless} {} {} {\itshape} {} {\bfseries} {.} { } {#1#3}
\theoremstyle {numless} 
\theoremstyle {definition} \newtheorem {df}   {Definition}
\theoremstyle {remark}  

\DeclareMathOperator {\ad}     {{ad}}

\DeclareMathOperator {\Aut}    {{Aut}}

\DeclareMathOperator {\Card}   {{Card}}
\DeclareMathOperator {\CH}     {{CH}}
\DeclareMathOperator {\Exp}    {{Exp}}
\DeclareMathOperator {\Fun}    {{Fun}}

\DeclareMathOperator {\Hom}    {{Hom}}
\DeclareMathOperator {\id}     {{id}}
\DeclareMathOperator {\im}     {{im}}

\DeclareMathOperator {\Log}    {{Log}}
\DeclareMathOperator {\Nil}    {{Nil_L}}
\DeclareMathOperator {\nilp}   {{\mathfrak{nilp}_L}}

\DeclareMathOperator {\pr}     {{pr}}
\DeclareMathOperator {\perv}   {{perv}}
\DeclareMathOperator {\rank}   {{rank}}

\DeclareMathOperator {\Rep}    {{Rep}}

\DeclareMathOperator {\Ext}    {{Ext}}
\DeclareMathOperator {\Loc}    {{Loc}}

\newcommand {\fa} {\mathfrak {a}}
\newcommand {\fb} {\mathfrak {b}}
\newcommand {\fc} {\mathfrak {c}}

\newcommand {\fg} {\mathfrak {g}}
\newcommand {\fh} {\mathfrak {h}}

\newcommand {\fr} {\mathfrak {r}}
\newcommand {\fp} {\mathfrak {p}}
\newcommand {\ft} {\mathfrak {t}}

\newcommand {\bQ} {\mathbb {Q}}

\newcommand {\bZ} {\mathbb {Z}}

\newcommand {\cE} {\mathcal {E}}
\newcommand {\cF} {\mathcal {F}}
\newcommand {\ccL}{\mathcal {L}}
\newcommand {\cM} {\mathcal {M}}

\newcommand {\cU} {\mathcal {U}}

\newcommand {\QZp} {\bQ_p/\bZ_p}
\newcommand {\ql}  {\overline{\bQ}_\ell}

\newcommand {\ld}   {\lambda}

\newcommand {\tensor}   {\otimes}

\newcommand {\xrarrow}[1]{\buildrel{#1}\over{\longrightarrow}}

\newcommand {\p}   	{\prime}

\newcommand {\wh}   	{\widehat}
\newcommand {\wt}  	{\widetilde}
\newcommand {\ov}  	{\overline}
\newcommand {\project} 	{\twoheadrightarrow}
\newcommand {\surject} 	{\twoheadrightarrow}
\newcommand {\inject}  	{\hookrightarrow}
\newcommand {\xarrow}[2] {\buildrel{#2}\over{#1}}
\newcommand {\Perf}     {\mathfrak{perf}}

\newcommand {\x}    {\chi}
\newcommand {\E}    {E}

\newcommand {\der}  {\mathcal{D}}

\newcommand {\bs}   {\backslash}
\newcommand {\wx}   {\wt{\chi}}

\title {Modular categories, orbit method and character sheaves on unipotent groups}
\author {Swarnendu Datta}
\date {}

\begin {document}
\maketitle

\begin {abstract}

Let $G$ be a unipotent group over a field of characteristic $p > 0$.
The theory of character sheaves on $G$ was initiated by 
V.~Drinfeld and developed jointly with D.~Boyarchenko. They also
introduced the notion of
$\mathbb{L}$-packets of 
character sheaves. Each
$\mathbb{L}$-packet can be described in terms of a modular
category. Now suppose that the nilpotence class of $G$ is less than $p$. Then the
$\mathbb{L}$-packets are in  bijection with the set $\mathfrak{g}^*/G$
of coadjoint orbits, where $\mathfrak{g}$ is the Lie ring scheme obtained from
$G$ using the Lazard correspondence and $\mathfrak{g}^*$ is the Serre dual of
$G$. If $\Omega$ is a coadjoint orbit, then the corresponding modular
category can be identified with the category of $G$-equivariant local
systems on $\Omega$. This in turn is equivalent to the category of
finite dimensional representations of a finite group. However, the
associativity, braiding and ribbon constraints are
nontrivial. Drinfeld gave a conjectural description of these
constraints in 2006. In this article, we prove the formula describing the ribbon structure when
$\dim(\Omega)$ is even.  
\end {abstract}
\tableofcontents

\section {Introduction}

Let $k$ be an algebraically closed field of characteristic $p > 0$ and
let $G$ be a unipotent group over $k$. Fix a prime $\ell$ different
from $p$ and let $\der(G)$ denote the bounded derived category of
constructible complexes of $\ell$-adic sheaves on $G$. Let $\der_G(G)$
be the category of $G$-equivariant objects of $\der(G)$, where $G$
acts on itself by conjugation.
The category
$\der_G(G)$ is monoidal with respect to convolution with compact
support. A (weak)
idempotent in  $\der_G(G)$ is an element $e
\in \der_G(G)$ such that $e * e \cong e$. A nonzero idempotent $e$ is
called {\em minimal} if for every idemponent $e^\p$, one has  $e *
e^\p \cong e$ or 0. The following fundamental
result is proved in \cite{foundations} (cf. Theorem 1.15): If $e \in
\der_G(G)$ is a minimal idempotent, then $e\der_G(G)$ is the
derived category of a modular category $\cM_e$ (a precise definition
of $\cM_e$ is given in Section 2.1). Two fundamental
problems in the theory of character sheaves over unipotent groups
are:\\ 

1. Give a concrete description of the
  set $\wh{G}$ of isomorphism classes of minimal
  idempotents in $\der_G(G)$.

2. For each $e \in \wh{G}$, describe the modular category $\cM_e$.\\

Partial solutions of the above problems are known in the general
situation (with arbitrary $G$). Essential to this appoach is a class
of special minimal idempotents called Heisenberg
  idempotents. These idempotents as well as the corresponding modular
  categories can be described explicitly. Moreover, one can show that
  the study of arbitrary idempotents can be 
reduced such idempotents. For an exposition of this
theory, we refer to \cite{foundations} and \cite{tanmay}.

Now assume that $G$  is connected and the nilpotence class of $G$ is
less than $p$. 
 In this case, a suitable modification of the orbit method applies and
 a complete answer 
to problem (1) is known. Let $\fg$ be the Lie ring scheme associated 
to $G$ via the Lazard correspondence and let $\fg^*$ be the Serre
dual of $\fg$. Using Fourier-Deligne transform, one can identify
$\wh{G}$  with
the set $\fg^*/G$ of coadjoint orbits (cf. Section 2.5).  
Let $e \in \der_G(G)$ be a minimal idempotent and let $\Omega \in
\fg^*$ be the coadjoint orbit corresponding to $e$. It can be shown
that $\cM_e$ is equivalent to the category of $G$-equivariant
local systems on $\Omega$. If $\x \in \Omega$ and $G_\x$ is the
stabilizer of $\x$ in $G$, then the later category is isomorphic to
$\Rep(\Gamma)$, where $\Gamma = \pi_0(G_\x)$. Thus one has a
equivalence  $\cM_e \simeq \Rep(\Gamma)$. This equivalence transforms
the monoidal functor on $\cM_e$ to usual tensor product on
$\Rep(\Gamma)$. However, the ribbon, braiding and associativity
constraints in $\cM_e$ induce 
nontrivial constraints in $\Rep(\Gamma)$. If $A = \ql[\Gamma]$, then
these constraints are given by certain elements in $A^\times,\, (A
\tensor A)^\times$ and $(A \tensor A \tensor A)^\times$. A conjectural
description of 
these elements were proposed by Drinfeld. The main result of
this paper is a proof of his formula for
the ribbon element in $A^\times$. The precise statement appears in
Theorem \ref{ribbon} in Section 3.3.\\

{\bfseries Acknowledgement}. I thank V. Drinfeld for sharing his
conjecture and for helpful discussions.

\section {The orbit method}
Let $k$ be an algebraically closed field of
characteristic $p > 0$ and 
let $G$ be a connected unipotent group over $k$. If the nilpotence
class of $G$ less than $p$, then one can use
the orbit method to study character sheaves on $G$. In this section,
we briefly recall the relevent ideas and constructions from
\cite{motivated}.

In the
remainder of this 
article we shall assume (unless otherwise stated) that:  
\begin {convention} The nilpotence class of $G$ is less
  than $p$ and $G$ is a perfect unipotent group {\em (cf. Section 1.9,
\cite{foundations})}.
\end {convention}
The first condition is necessary for the orbit
method to work. The second 
condition is imposed to avoid technical difficulties involving Serre
duality (of perfect commutative unipotent groups - cf. Section 2.3
below). However, the results of this paper are independent of this
condition, as is explained in \cite {foundations}.

Fix a prime
$\ell$ different from $p$. For a scheme $X$ over $k$, let $\der(X)$
denote the bounded derived category of constructible complexes of
$\ell$-adic sheaves on $X$. If $X$ is equipped with a $G$ action, then
let $\der_G(X)$ denote the category of $G$-equivariant objects of
$\der(X)$ (cf. Appendix E, \cite{motivated}).

\subsection {The category $\der_G(G)$ and character sheaves on $G$}
In this subsection, we give the definition of the category $\cM_e$
mentioned in the introduction. The assumption on the nilpotent class
of $G$ is not used.

Let $\der_G(G)$ be the category of $G$-equivariant objects of
$\der(G)$, where $G$ acts on itself by conjugation. 
Let $M$, $N$ be objects on $\der_G(G)$. The convolution of $M$ and $N$
is the object $M*N$ of $\der_G(G)$ defined by
$$
M*N = \mu_!(\pr_1^*(M) \otimes \pr_2^*(N))
$$
where $\mu: G \times G \rightarrow G$ is the multiplication map and
$\pr_1$, $\pr_2: G \times G$ are projections onto the first and second
factors. The category $\der_G(G)$ is a monoidal category with respect
to the convolution bifunctor. Moreover, one can define natural
isomorphims
$$
\beta_{M,N}: M*N  \rightarrow N*M
$$
which equips $\der_G(G)$ with the structure of a braided
category. These isomorphisms are constructed using the $G$
equivariance of $N$ (cf. Section 5.5, \cite{motivated}). Further,
there is a canonical  automorphism $\theta$ of the identity
functor of $\der_G(G)$. If $M \in \der_G(G)$, then the automorphism
$\theta_M:M \rightarrow M$ can be described as follows. For every $g$,
$x \in G$, one has isomorphisms $\phi_{g,x}: M_g \rightarrow
M_{xgx^{-1}}$. Setting $g = x$ gives an automorphism of $M_g$, which
is none other than $(\theta_M)_g$. More precisely, if $c: G \times G
\rightarrow G$ is defined by $c(x,g) = xgx^{-1}$, then the $G$
equivariant structure on $M$ yeilds an isomorphism $\pr_2^*
\rightarrow c^*M$. Pulling back this isomorphism by the diagonal map 
$\Delta: G \rightarrow G \times G$ gives $\theta_M$.  One has the
following relation for all $M$, $N \in \der_G(G)$:
$$
\theta_{M*N} = \beta_{N,M} \circ \beta_{M,N} \circ (\theta_M *
\theta_N)
$$
This can be checked directly. An abstract proof is given in
Proposition A.46, \cite{foundations}.

Let $e$ be a minimal idempotent of $\der_G(G)$. The modular category
$\cM_e$ mentioned in the introduction is defined as follows. Let
$e\der_G(G)$ be the full subcategory of $\der_G(G)$ consisting of
objects of the form $e*M$ where $M \in \der_G(G)$. Let $\cM_e^{\perv} \subset
e\der_G(G)$ denote the full subcategory of perverse objects in
$e\der_G(G)$. A {\em character sheaf} on $G$ is defined to be an
indecomposible object of $\cM_e^{\perv}$ for some minimal idempotent
$e$. It is shown in \cite{foundations} that there exists a
(unique) integer $n_e \geq 0$ such that $e[-n_e] \in \cM_e^{\perv}$. One now
defines $\cM_e$ to be the subcategory $\cM_e^{\perv}[n_e]$ of
$\der_G(G)$. It is closed under convolution and is a modular category,
with braiding and ribbon structure given by $\beta$ and $\theta$
(cf. Theorem 1.15, \cite{foundations}).  

\subsection {Lazard correspondence}

A {\em Lie ring} is an additive group $\fp$ equipped a biadditive map
$[\cdot ,\cdot]: \fp \times \fp \rightarrow \fp$ satisfying the Jacobi
identity and the identity $[x,x] = 0$. Let $\nilp$ be the category of
all nilpotent Lie rings (i.e.,
Lie algebras over $\mathbb{Z}$) $\fp$ satisfying the following
condition:\\ 

If $n$ is the nilpotence class of $\fp$, then the map $\fp
\xrightarrow {\times n!} \fp$ of multiplication by $n!$ is
invertible. \hfill (L)\\ 

\noindent One can define a multiplication operation $\fp \times \fp
\rightarrow \fp$, $(x,y) \mapsto xy$ by using the Campbell-Hausdorff
formula 
$$
xy = \log(e^x\cdot e^y) = x + y + \frac{1}{2}[x,y] +
\frac{1}{12}([x,[x,y]] + [y,[y,x]]) + \cdots 
$$
If $\CH_i$ is homogeneous component of degree $i$ is the above series,
then one knows that $\CH_i$ is a Lie polynomial over
$\mathbb{Z}[1/i!]$ (cf. Section IV.8 of \cite{serre_lie}). The
condition (L) ensures that the above formula is well defined. Let
$\Exp(\fp)$ be the group whose underlying set is same as that of $\fp$
and whose product operation is given by the above formula. The inverse
of $x$ is $-x$ and the identity element is $0 \in \fp$. If
$\varphi:\fp_1 \rightarrow \fp_2$ is a homomorphism of Lie rings in
$\nilp$, then it is immediate that $\varphi$, viewed as a map
$\Exp(\fp_1) \rightarrow \Exp(\fp_2)$, is a group homomorphism. In
other words, $\Exp$ is a functor from $\nilp$ to the category of
groups. Let $\Nil$ be the category of all nilpotent groups $\Gamma$
such that if $n$ is the nilpotence class of $\Gamma$, then the map
$\Gamma \rightarrow \Gamma$, $g \mapsto g^{n!}$ is invertible.  

\begin {thm*}[Lazard]
If $\fp \in \nilp$ then $\Exp(\fp) \in \Nil$. The functor $\Exp: \nilp
\rightarrow \Nil$ is an isomorphism of categories. 
\end {thm*}

The inverse of the functor $\Exp$ is denoted by $\Log$.
A proof of the theorem can be found in \cite{lazard} or
\cite{khukhro}. The difficult 
part is to show that one can recover $\fp$ from $\Exp(\fp)$, i.e, the
addition and Lie bracket operations can be expressed in terms of the
product operation. For example, in nilpotence class 2, it can be
checked that $[x,y] = (x,y)$,\, $x + y = xy(x,y)^{-1/2}$ where $(x,y)$
denotes $xyx^{-1}y^{-1}$. 

For $x \in \Exp(\fp)$, the conjugation
action of $x$ on  $\Exp(\fp)$
gives an automorphism of $\Exp(\fp)$. By the theorem of Lazard, it is
also an automorphism of the Lie ring $\fp$. This is easy to check
directly using the following formula:
\begin {lem}
\label {conjugateformula}
$xyx^{-1} = e^{\ad(x)}(y) = y + [x,y] + \frac{1}{2}[x,[x,y]] + \frac{1}{6}[x,[x,[x,y]]]+\cdots$
\end {lem}
\begin {proof}
We have $e^{xyx^{-1}} = e^x\cdot e^y \cdot e^{-x} = 1 + e^x\cdot y
\cdot e^{-x} + \frac{1}{2}e^x\cdot y^2 \cdot e^{-x} + \cdots =
e^{e^x\cdot y \cdot e^{-x}}$. An easy computation shows that
  $e^x\cdot y \cdot e^{-x}
= x + [x,y] + \frac{1}{2}[x,[x,y]] + \cdots$, qed.
\end {proof}
One can use the Lazard correspondence to produce a Lie ring scheme
(i.e., a Lie ring object in the category of schemes) from
$G$. Indeed, for every integer $n$ with $1 \leq n \leq p-1$, the
multiplication by $n$ map from $G$ to itself is invertible (this
follows by using a filtration of $G$ such that the successive
quotients are isomorphic to the additive group). Hence for every
scheme $S$ over $k$, the group $G(S)$ belongs to $\Nil$. Hence the
set $G(S)$ can be given the structure of a Lie ring. Therefore one
obtains the structure of a Lie ring scheme on $G$. This Lie ring
scheme is denoted by $\Log(G)$. Note that $G$ and $\Log(G)$ are {\em
  equal} as schemes.

\subsection {Serre duality and mulplicative local systems}
We refer to Appendix F, \cite{motivated} for a more detailed
exposition. 
Let $\Perf_k$ denote the category of perfect schemes over $k$ and let
$A$ be a commutative, unipotent, connected group scheme in
$\Perf_k$. One considers the following functor, from the category
$\Perf_k$ to abelian groups: 
$$
S \mapsto \Ext^1(A \times S, \QZp)
$$
where $\Ext^1$ is computed in the category of commutative group
schemes over $S$ and $\QZp$ is viewed as a discrete group scheme over
$S$. One knows that the functor defined
above is  representable by a group
scheme $A^*$ in $\Perf_k$ (cf. \cite{begueri}, Prop. 1.2.1). The
group $A^*$ is called the \emph{Serre dual} of $A$. One has the
following properties (\emph{loc. cit.}): 

i) $A^*$ is a connected, commutative, unipotent group scheme isogenous
to $A$. 

ii) The canonical homomorphism of $A$ onto $A^{**}$ is an isomorphism.  

iii) If $0 \rightarrow A^\p \rightarrow A \rightarrow A^{\p\p}
\rightarrow 0$ is an exact sequence of connected, commutative,
unipotent group schemes in $\Perf_k$, then so is $0 \rightarrow
A^{\p\p*} \rightarrow A^* \rightarrow A^{\p*} \rightarrow 0$.

iv) If $f: A \rightarrow A'$ is an isogeny, then so is $f^*: A^{\p*}
\rightarrow A^*$.

One also has the following alternative intepretation of the Serre
dual. Let $H$ be a group scheme over $k$ and let $\ccL$ be a
$\ql$-local system on  $H$. If $\mu:H \times H \rightarrow H$ is the
multiplication map and $\pr_1$, $\pr_2: H\times H \rightarrow H$ are
the projection maps, then $\ccL$ is called a {\em multiplicative} local
system if $\ccL$ is nonzero and $\mu^*\ccL \cong \pr_1^*\ccL \otimes
\pr_2^*\ccL$.  
 Fix an injective homomorphism (of groups):
$$
\psi: \QZp \inject \ql^\times
$$
An element $\x \in A^*$ is a (central) extension of $A$ by
$\QZp$. 
 Using the embedding $\psi:\QZp \inject \ql^\times$, one can view $\x$
 as a multiplicative $\ql$-local system on $A$. Conversely, it can be
 proved that every multiplicative $\ql$-local system on $A$ comes from
 an element of $A^*$ via $\psi$ (cf. Lemma 7.3 in \cite{mitya}). 

\subsection {Fourier-Deligne transform}
Let $\fg$ denote the Lie ring scheme $\Log(G)$ (cf. Section 2.2)
obtained from $G$ using the Lazard correspondence. In particular,
$\fg$ is a commutative group scheme and thus its Serre dual $\fg^*$ is 
well defined. 
Let $\cU$ be the universal $\QZp$ torsor on $\fg \times \fg^*$. If $g
\in \fg$ and $\x \in \fg^*$ then the stalk of $\cU$ at $(g,\x)$ is the
stalk of $\x$ at $g$. Let $\cE$ be the local system on $\fg \times
\fg^*$ obtained from $\cU$ using the homomorphism $\psi$. The
Fourier-Deligne transform is the functor 
$$
\cF: \der(\fg) \rightarrow \der(\fg^*)
$$
defined by
$$
\cF(M) = \pr'_!(\pr^*(M) \otimes \cE)
$$
where $\pr$, $\pr'$ are projections: $$\fg \xleftarrow{\pr'} \fg \times
\fg^* \xrightarrow {\pr} \fg$$
One also has the functor $\cF': \der(\fg^*) \rightarrow \der(\fg)$
defined by 
$\cF'(N) = \pr_!(\pr'^*(N)\otimes \cE)$. Let $d = \dim(\fg)$. For each
$M \in \der(\fg)$, there exist a natural isomorphism 
$$
(\cF'\circ \cF) (M) \cong (-1)^*M[-2d](-d)
$$
where $-1: \fg \rightarrow \fg$ denotes the inverse map. This is
proved in \cite{saibi} (cf. Th\'eor\'eme 2.2.4.1). It follows that
$\cF$ is an equivalence of categories and a quasi-inverse functor is
given by $N \mapsto (-1)^*\cF'(N)[2d](d)$.
Further, $\cF$ transforms convolution in $\der(\fg)$ (with respect to
the additive group structure in $\fg$) to tensor product in
$\der(\fg^*)$. That is, for all $M$, $N \in \der(\fg)$, there exist
natural bifunctorial isomorphisms
$$
\cF(M*N) \cong \cF(M)\otimes\cF(N)
$$
inducing an equivalence of the symmetric monoidal categories
$(\der(\fg),*)$ and $(\der(\fg^*),\otimes)$. This is proved in
Proposition H.11, \cite{motivated}. 

\subsection {The bijection between $\wh{G}$ and $\fg^*/G$}
Consider the action of $G$ on itself by conjugation. Since $G$ and $\fg$
are equal as schemes, this induces an action of $G$ on $\fg$. One also
has the corresponding contragradient action of $G$ on $\fg^*$. It is
proved in \cite{motivated} (cf. Appendix H) that $\cF$ can be lifted
to an equivalence 
$$
\cF : \der_G(\fg) \rightarrow \der_G(\fg^*)
$$
of the monoidal categories $(\der_G(\fg),*)$ and
$(\der_G(\fg^*),\otimes)$.

In the following, we shall identify $\der_G(G)$ and $\der_G(\fg)$. For
$M$, $N \in \der_G(G)$, let $M *_{\text{a}} N$ denote that
additive convolution with respect to $\fg$ and let $M *_{\text{m}}
N$ denote the multiplicative convolution with respect to $G$. The
important observation is that there exist bifunctorial isomorphisms 
\begin {equation}
M *_{\text{a}} N \cong M*_{\text{m}} N
\tag{$*$}
\end {equation}
for all $M$, $N \in \der_G(G)$. These isomorphims are
not canonically defined and depend on the choice of certain Lie
polynomials 
(cf. Proposition 5.11, {\em loc. cit.}). However, it follows that an
element $e \in (\der_G(G), *_{\text{m}})$ is a minimal idemponent if
and only if it is 
so in $(\der_G(\fg), *_{\text{a}})$. Via the Fourier-Deligne transform,
such idemponents correspond to minimal idempotents in
$(\der_G(\fg^*),\otimes)$. Let $\Omega$ be a $G$ orbit in
$\fg^*$. Since $G$ is unipotent, one knows that $\Omega$ is closed in
$\fg^*$. Let 
$i_\Omega:\Omega \inject \fg^*$ be the inclusion map and let
$(\ql)_\Omega$ be the constant sheaf on $\Omega$ with stalk $\ql$. Then
$e_\Omega = (i_\Omega)_*((\ql)_\Omega)$ is a minimal idempotent of
$\der_G(\fg^*)$. Moreover, all such idempotents are of the form
$e_\Omega$ for a suitable orbit $\Omega \subset \fg^*$. Therefore one
obtains a bijection between $\wh{G}$ and $\fg^*/G$. 

Let $e$ be the
minimal idempotent of $\der_G(G)$ corresponding to an orbit
$\Omega$. One has $\cF(e) \cong e_\Omega$. The Fourier-Deligne transform
gives an equivalence between $e\der_G(G)$ and $e_\Omega\der_G(\fg^*)$,
which can be identified with $\der_G(\Omega)$. This
equivalence transforms
convolution in $e\der_G(G)$ to usual tensor product in
$\der_G(\Omega)$. However, for general $G$, the braiding and
associativity constraints in $e\der_G(G)$ are transformed into
nontrivial constraints in $\der_G(\Omega)$. This is because the
isomorphims in $(*)$ above are not compatible (in general) with the
braiding 
and associativity constraints in $\der_G(G)$ and $\der_G(\fg)$.
Let $\Loc_G(\Omega)$ denote the full
subcategory of $\der_G(\Omega)$ consisting of $G$-equivariant local
systems on $\Omega$. 

\begin {lem}
\label {equivalence}
The Fourier-Deligne transform induces an equivalence between $\cM_e$
and $\Loc_G(\Omega)$.
\end {lem}

\begin {proof}
The category of the perverse objects in $\der_G(\Omega)$ is
$\Loc_G(\Omega)[\dim \Omega]$. Let $d = \dim G$. As $\cF[d]$ takes
perverse sheaves to 
perverse sheaves (cf Appendix H, \cite{motivated}), it follows that
$\cM_e^{\perv}$ is  equivalent to $\Loc_G(\Omega)[\dim \Omega -
d]$. Therefore $\cM_e^{\perv}[d - \dim \Omega]$ contains $e$ and thus $\cM_e =
\cM_e^{\perv}[d - \dim \Omega]$, which is equivalent to
$\Loc_G(\Omega)$, qed.
\end {proof}

\section {Statement of the main result}

\subsection {Skew-symmetric biextensions}
We recall several facts about biextensions
from Appendix A of \cite{mitya}. For the definition of biextensions in
general, the reader is refered to {\em loc. cit.}. 
Let $\fa$, $\fb$ be perfect, commutative, connected, unipotent groups and
let $E$ be a biextension of $\fa \times \fb$ by $\QZp$. In particular,
$E$ is a 
$\QZp$ torsor over $\fa \times \fb$. For each
$y \in \fb$, the restriction $E|_{\fa \times y}$ is an extension
  of $\fa$ by $\QZp$, and thus is an element of $\fa^*$. The map $y
  \mapsto E|_{\fa \times y}$ defines an additive homomorphism $f: \fa
  \rightarrow \fb^*$. One can similarly define a map from $\fb$ to
  $\fa^*$, which is 
 none other than the dual $f^*$ of $f$. The restriction of $E$ to
 $\ker f \times \fb$ and $\fa \times \ker f^*$ are trivial and thus one
 obtains unique trivializations $\ld: \ker f \times \fb \rightarrow E$
 and $\rho: \fa \times \ker f^* \rightarrow E$ such that $\ld(0,0) = 0 =
 \rho(0,0)$. Define a map $\wt{B}: \ker f \times \ker f^* \rightarrow
 \QZp$ by $\wt{B}(x,y) = \ld(x,y) - \rho(x,y)$. Since $\QZp$ is
 discrete, this descends to a (biadditive) pairing
$$
B: \pi_0(\ker f) \times \pi_0(\ker f^*) \rightarrow \QZp
$$
It is proved in Proposition A.19, {\em loc. cit.} that $B$ is
nondegenerate. Now suppose that $E$ is a {\em skew-symmetric}
biextension of $\fa \times \fa$, i.e., the restriction of $E$ to the
diagonal $\Delta(\fa) \subset \fa \times \fa$ is trivial. In this case
one has $f = -f^*$ and thus $\ker f = \ker f^*$. Let
$\alpha: \Delta(\fa) \rightarrow E$ be the trivialization of
$E|_{\Delta(\fa)}$ such that $\alpha(0) = 0$. Define $\wt{q}:
\ker f \rightarrow \QZp$ by $\wt{q}(x) = \ld(x,x) - \alpha(x,x)$. Then
$\wt{q}$ descends to a map $q: \pi_0(\ker f) \rightarrow \QZp$. One
has the following result (cf. Lemma A.26, {\em loc. cit.}):

1) $q(nx) = n^2q(x)$ for all $n \in \mathbb{Z}$.

2) $B(x,y) = q(x+y) - q(x) - q(y)$.

\noindent In particular, $q$ is a nondegenerate quadratic form on $\pi_0(\ker
f)$. 

\subsection {The quadratic form $q$}
Let $\Omega$ be a $G$-orbit in $\fg^*$. Fix an element $\x \in
\Omega$. Let $G_\x$ be the stabilizer of $\x$ in $G$. Let $\Gamma =
\pi_0(G_\chi)$ and let $\fp = \Log(\Gamma)$. Using the construction of
the 
previous section, we shall now define a quadratic form $q: \fp
\rightarrow \QZp$ on $\fp$. This quadratic form is used in Drinfeld's
formula for the ribbon element in $\ql[\Gamma]$ (cf. Section 3.3).

Let $\E$ be the biextension of $\fg \times \fg$ obtained by
pulling back $\x$ via the biadditive commutator map 
$
\fg \times \fg \rightarrow \fg
$
given by $(u,v) \mapsto [u,v]$. The restriction of $\E$ to the
diagonal $\Delta\fg \subset \fg \times \fg$ is trivial, whence $\E$ is
a skew-symmetric biextension. If $f: \fg \rightarrow \fg^*$ is the
homomorphism corresponding to $E$, then the kernel of $E$ is defined
to be $\ker f$. Thus $x \in \ker E$ if and only if $E|_{x \times \fg}$
is trivial. Let $\fg_\x = \Log(G_\chi)$.

\begin {lem}
\label {kernel}
$\fg_\x$ is the kernel of $\E$.
\end {lem}

If $x \in G$, then we denote the element $x(\chi) \in \fg^*$ by
${}^x\chi$.

\begin {proof}
In the following, we identify $G$ and $\fg$ (resp. $G_\x$ and
$\fg_\x$) as schemes. Let $\fc$ denote the kernel of $E$. It suffices
to show that $\fc \subset \fg_\x$ and $\fg_\x \subset \fc$. 

(i) $\fc \subset \fg_\x$. Let $x \in \fc$. Then the restriction of
$\E$ to $x \times \fg$ is trivial. Equivalently, if $\phi_x: \fg
\rightarrow \fg$ is the map defined by $\phi_x(y) = [x,y]$, then
$\phi_x^*(\x)$ is trivial. We want to show that $x \in \fg_\x$, i.e,
$\x = {}^x\x$ in $\fg^*$. The extension ${}^x\x$ is obtained by
pulling back $\x$ by the 
automorphism $\varphi_x:\fg \rightarrow \fg$ defined by $\varphi_x(y)
= x^{-1}yx$. Thus $\x - {}^x\x = (\id -\, \varphi_x)^*(\x)$.
Further, using Lemma \ref{conjugateformula} one obtains:
\begin {align*}
(\id -\, \varphi_x)(y) = y - x^{-1}yx = y - e^{\ad (-x)}(y) &= [x,y] -
\tfrac{1}{2}[x,[x,y]] + \tfrac{1}{6}[x,[x,[x,y]]] + \cdots\\ 
&= [x,\mu(y)]
\end {align*}
where $\mu: \fg \rightarrow \fg$ is the additive map defined by
$\mu(y) = y - \tfrac{1}{2}[x,y] + \tfrac{1}{6}[x,[x,y]] - \cdots$. It
follows that $\id -\, \varphi_x$ can be written as the composition: 
$
\fg \xrightarrow{\mu} \fg \xrightarrow{\phi_x} \fg
$.
By hypothesis, $\phi_x^*(\x)$ is trivial, whence so is $(\id
-\,\varphi_x)^*(\x)$. Therefore $x \in \fg_\x$. 

(ii) $\fg_\x \subset \fc$. Let $x \in \fg_\x$. Then $(\id -
\varphi_x)^*(\chi) = 0$. To show that
$x$ is in the kernel of $\E$, it suffices to find (additive) $\lambda:
\fg \rightarrow \fg$ such that $\phi_x = (\id -\, \varphi_x) \circ
\lambda$, i.e., 
\begin {equation}
\label {eqn}
[x,y] = [x,\lambda(y)] - \tfrac{1}{2}[x,[x,\lambda(y)]] + \tfrac{1}{6}[x,[x,[x,\lambda(y)]]] + \cdots
\tag {$*$}
\end {equation}
Set
$$
\lambda(y) = y + \tfrac{1}{2}[x,y] + \tfrac{1}{12}[x,[x,y]] - \tfrac{1}{720}[x,[x,[x,[x,y]]]] + \cdots
$$
where the coefficients $B_n$ of $(\ad x)^n(y)$ satisfy the equation:
$$
B_0 + B_1t + B_2t^2 + \cdots = \frac{t}{1 - e^{-t}}
$$
It is easily checked that $\lambda(y)$ is well
defined (i.e., the denominator of $B_n$ is coprime to the $p$ whenever
$n+1 < p$) and that it satisfies (\ref{eqn}). This completes the
proof. 
\end {proof}

Let $\fp = \pi_0(\fg_\x)$. Note that $\fp = \Log(\Gamma)$. Applying
the construction of Section 3.1 to $E$, one obtains a quadratic form:
$$
q: \fp \rightarrow \QZp
$$

\begin {lem}
\label {invariant}
$q$ is invariant under the conjugation action of $\Gamma$ on $\fp$.
\end {lem}

\begin {proof}
Let $\wt{q}: \fg_\x \rightarrow \QZp$ be the function obtained by
pulling back $q$ via the projection $\fg_\x \project \fp$. It suffices
to show that $\wt{q}$ is invariant under the action of $G_\x$. 
Let $x$ be an element of $G$ and let ${}^x\E$ be the biextension of
$\fg \times \fg$ obtained from ${}^x\x \in \fg^*$ as above. If
$\varphi_x$ is as in the proof of Lemma \ref{kernel}, then ${}^xE =
(\varphi_x \times \varphi_x)^*(E)$. Thus
$
\ker ({}^x\E) = \varphi^*_x (\ker E) = x\fg_\x x^{-1}
$.
The corresponding function $\wt{q}_x: \ker ({}^x\E) \rightarrow \QZp$
is $\wt{q} \circ \varphi_x$. Now assume that $x
\in G_\x$. Then ${}^x\x = \x$, whence ${}^x\E \cong \E$. Thus
$\wt{q}_x = \wt{q}$, which proved the desired assertion.
\end {proof}

\subsection {Drinfeld's formula}
The category
$\Loc_G(\Omega)$ of $G$-equivariant local systems on $\Omega$ is
equivalent to the category $\Rep(\Gamma)$ of finite dimensional
representation of $\Gamma$ over $\ql$ (the equivalence sends a local
system $\ccL$ to its stalk $\ccL_\chi$ at $\chi$). It
follows from Lemma \ref{equivalence} that $\cM_e$ is equivalent to
$\Rep(\Gamma)$. The ribbon structure of $\cM_e$ given by $\theta$
induces an automorphism of the identity functor of
$\Rep(\Gamma)$. This automorphism can be interpreted as an invertible
element of the center of $\ql[\Gamma]$. A precise formula (due to
Drinfeld) for this 
element is given in Theorem \ref{ribbon} below. 

Let $\wt{q}$ be the composition:
$$\fp \xrarrow{q} \QZp \xarrow{\hookrightarrow}{\psi} \ql^\times$$
Thus $\wt{q}$ is a quadratic form on $\fp$
with values in $\ql^\times$. Let $\fp^* = \Hom(\fp,\ql^\times)$ be the
group of characters of $\fp$. Since $\wt{q}$ is nondegenerate, the
bimultiplicative pairing $\fp \times \fp \rightarrow \ql^\times$
defined by $\wt{q}$ 
gives an isomorphism of the abelian groups $\fp$ and $\fp^*$. This
isomorphism transforms $\wt{q}$ into a quadratic form on $\fp^*$,
which we denote by $\wt{q}'$. It
easily follows from Lemma \ref{invariant} that $\wt{q}'$
is invariant under the action of $\Gamma$ on $\fp^*$. Consider the 
Fourier transform
$$
\ql[\fp] \rightarrow \Fun(\fp^*,\ql)
$$
which sends $a \in \fp$ to the function $\phi \mapsto
\phi(a)^{-1}$. Let
$\wh{q}$ be the preimage of $\wt{q}'$ under this map. Since
this map is $\Gamma$ invariant, it follows that $\wh{q}$ is also
$\Gamma$ invariant. Using the identity map $\Gamma
\rightarrow \fp$, one can interpret $\wh{q}$ as an element of the
center of $\ql[\Gamma]$.

\begin {thm}
\label {ribbon}
Let $\theta'$ be the automorphism of the identity functor of
$\Rep(\Gamma)$ corresponding to the ribbon automorphism $\theta$ under
the equivalence $\cM_e \simeq \Rep(\Gamma)$. Then $\theta'$ is
multiplication by $\wh{q}$.
\end {thm}

Conjecturally, the braiding and associativity constraints in
$\Rep(\Gamma)$ coming from $\cM_e$ is also completely determined by
$\wt{q}$. We now give an explicit formula for $\wh{q}$. Let $|\fp|$
denote the 
cardinality of $\fp$ and let $G(\fp,\wt{q}) = \sum_{a \in \fp}
\wt{q}(a)$. Note 
that $G(\fp,\wt{q})$ is the {\em Gauss sum} of the {\em metric group}
$(\fp,\wt{q})$ (cf. Section  6 of \cite{dgno}).

\begin {lem}
\label {formulaforqhat}
$\wh{q} = \dfrac{G(\fp,\wt{q})}{|\fp|} \displaystyle\sum_{a \in \fp}
\wt{q}(a)^{-1}a$ 
\end {lem}

\begin {proof}
For $\phi \in \fp^*$, let $1_\phi \in \Fun(\fp^*,\ql)$ be the
characteristic function of $\phi$. The inverse image of $1_\phi$ under
the Fourier transform is $1/|\fp| \sum_{a \in \fp} \phi(a)a$. Let
$\phi_a \in \fp^*$ be the character corresponding to $a \in
\fp$ (under the isomorphism $\fp \rightarrow \fp^*$ defined using
$\wt{q}$). Then  $\wt{q}' = \sum_{b \in \fp} \wt{q}(b)1_{\phi_b}$. If
$\wt{B}:\fp \times \fp \rightarrow \ql^\times$ is the bimultiplicative
pairing defined by $\wt{q}$, then one has $\phi_b(a) = \wt{B}(a,b)$. Thus 
$$
\wh{q} = \frac{1}{|\fp|} \sum_{b \in \fp}\Bigl( \wt{q}(b)\sum_{a\in
    \fp}\wt{B}(a,b)a\Bigr) 
$$
As $\wt{q}(b)\wt{B}(a,b) = \wt{q}(a)^{-1}\wt{q}(a+b)$, it follows that the
coefficient of 
$a$ in the above sum is $\wt{q}(a)^{-1}G(\fp,\wt{q})/|\fp|$, qed.
\end {proof}

\subsection {A reformulation of Theorem \ref{ribbon}}
Let $\ccL$ be the $G$-equivariant local system on $\Omega$
corresponding to the regular representation of $\Gamma$. To prove
Theorem \ref{ribbon}, it suffices to show that $\theta'$ acts by multiplication
by $\wh{q}$ on the stalk $\ccL_\chi$. Let $T$ be an element of $\cM_e$
corresponding to $\ccL$. A formula for $T$ is given in Lemma
\ref{inverseofl} using the inverse Fourier-Deligne transform. To
understand $\theta'_{\ccL}: \ccL_\x \rightarrow \ccL_\x$, we can use
the fact that it is obtained
by applying $\cF$ to $\theta_T: T \rightarrow T$. This approach yeilds an
explicit interpretation of $\theta'_\ccL$, described in Proposition
\ref{reformulation} below.\\

{\bfseries Notation}. If $\varphi$ is an element of $\fg^*$, then let
$\wt{\varphi}$  denote the (multiplicative) local system on $\fg$
obtained from $\varphi$ using $\psi:\QZp \inject \ql^\times$.\\

Let $$w: \fg \times G \rightarrow \fg$$ be the map $(x,y) \mapsto x -
y^{-1}xy$ and let $W = w^*(\wt{\chi})$. Then:

A) The map $w$ is invariant under the automorphism $(x,y) \mapsto
(x,xy)$ of $\fg \times G$.

B)  Consider the action of $G_\chi$ on $\fg \times G$ such that $g
\in G_\x$ takes $(x,y)$ to $(gxg^{-1},gy)$. Then $W$ has a natural
$G_\x$-equivariant structure. Indeed,
let $v: \fg \times G \rightarrow \fg$ be the map which
takes $(x,y)$ to $-y^{-1}xy$. Since $\wt{\chi}$ is multiplicative, one
has $W = \pr_1^*(\wt{\x})\otimes v^*(\wt{\x})$. Let $g \in G_\x$ and
let $\sigma_g$ be the corresponding automorphism of $\fg \times
G$. Then $$\sigma_g^*(W) = \pr_1^*(\wt{{}^{g^{-1}}\x})\otimes
v^*(\wt{\x})$$ This is isomorphic to $W$ since ${}^{g^{-1}}\x = \x$.

Let $V = H^{2i}_c(W)$, where $i = \dim G +
\dim G_\chi$ and $H^\bullet_c$ denotes $\ell$-adic cohomology with compact
support. It follows from (B) that $V$ is a $\Gamma$ module. The
automorphism in (A) induces an automorphism $\eta:V \rightarrow V$ of
$V$.

\begin {prop}
\label {reformulation}
There exists an isomorphsim $\ccL_\chi \cong V$ of $\Gamma$ modules
which maps $\theta'_\ccL$ to $\eta$.
\end {prop}

The proof of the proposition is given in Section 3.6. It follows that
Theorem \ref{ribbon} is equivalent to showing that $\eta$ acts by multiplication
by $\wh{q}$. This is proved by explicit computation in Sections 5, 6.

\subsection {A formula for $T$ and the automorphism $\theta_T$}
 We use the following notation: if $\mu:X \rightarrow Y$ is a
morphism, then $\mu_!(\ql)$ is $\mu_!((\ql)_X)$, where $(\ql)_X$ is
the constant sheaf on $X$ with stalk $\ql$. Let $d = \dim G$ and $n =  
\dim G_\chi$.  
\begin {lem} Let $\pi: G \rightarrow G/G_\chi = \Omega$ be the
projection map. Then $\pi_!(\ql) = \ccL[-2n](-n)$.
\end {lem}

\begin {proof}
The morphism $\pi$ is the composition of the projection maps:
$
G \xrightarrow{\pi_1} G/G_\chi^\circ \xrightarrow{\pi_2} G/G_\chi
$.
Note that $(\pi_1)_!(\ql) = \ql[-2n](-n)$. This follows since
$G_\x^\circ$ is isomorphic to the affine space of dimension $n$ and $G$  
is a $G_\x^\circ$-torsor over $G/G_\x^\circ$. Further, $\pi_2$ is a Galois
cover with Galois group $\Gamma$, whence $(\pi_2)_!(\ql) \cong
\ccL$. This completes the proof. 
\end {proof}

Let $\ov{\ccL}$ be the
sheaf on $\fg^*$ obtained by extending $\ccL$ by zero. Let $\cF^{-1}$
be a quasi-inverse of $\cF$ which sends $N$ to
$(-1)^*\cF'(N)[2d](d)$ (cf. Section 2.4). Then one has $T \cong
\cF^{-1}(\ov{\ccL})$.

\begin {lem}
\label {inverseofl}
Let $v: \fg \times G \rightarrow \fg$ be the map $(x,y) \mapsto
-y^{-1}xy$ and let $V = v^*(\wt{\chi})$. If $\pr_1: \fg \times G
\rightarrow \fg$ is the first projection, then $$T \cong
(\pr_1)_!(V)[2n+2d](n+d)$$
\end {lem}

\begin {proof}
Let $u:G \rightarrow \fg^*$ be the composition
$$
G \xrightarrow {\pi} \Omega \inject \fg^*
$$
One has $u_!(\ql) = \ov{\ccL}[-2n](-n)$. Consider the commutative
diagram:
$$
\xymatrix{
& \fg \times G \ar[r]^-{\pr_2} \ar[d]^-{\wt{u}} & G \ar[d]^-{u}\\
\fg & \fg \times \fg^* \ar[l]_-{\pr} \ar[r]^-{\pr'} & \fg^*
}
$$
where $\wt{u} = \id_\fg \times u$. The square is cartesian and
using the projection formula, it follows that
\begin {align*}
\pr'^*(\ov{\ccL})\otimes \cE &= (\pr'^*u_!(\ql) \otimes \cE)[2n](n)\\ 
& \cong (\wt{u}_!(\ql)\otimes \cE) [2n](n)\\
& \cong (\wt{u}_!\wt{u}^*(\cE))[2n](n)
\end {align*}
As $\pr_1 = \pr \circ \wt{u}$, one has
$$
\cF'(\ov{\ccL}) = \pr_!(\pr'^*(\ov{\ccL})\times \cE) \cong
((\pr_1)_!\wt{u}^*(\cE))[2n](n)
$$
 Let $\iota: \fg \times G
\rightarrow \fg \times G$ be defined by $(x,y) \mapsto (-x,y)$. Then 
\begin {align*}
T \cong \cF^{-1}(\ov{\ccL}) &\cong
(-1)^*(\pr_1)_!\wt{u}^*(\cE)[2n+2d](n+d)\\ 
&\cong (\pr_1)_!(\iota^*\wt{u}^*(\cE))[2n+2d](n+d)
\end {align*}
Now the stalk of $\wt{u}^*(\cE)$ at $(x,y)$ is $\cE_{x,u(y)}$. As
$u(y) = {}^y\x$, the stalk is $\wt{({}^y\chi)}_x =
\wt{\x}_{y^{-1}xy}$. Thus the stalk of 
$\iota^*\wt{u}^*(\cE)$ at $(x,y)$ is $\wt{\x}_{y^{-1}(-x)y} = 
  \wt{\chi}_{-y^{-1}xy}$. Therefore $\iota^*\wt{u}^*(\cE) = V$, qed.
\end {proof}
To understand the automorphism $\theta_T:T
\rightarrow T$, note that

1) The diagram appearing in the proof above is $G$-equivariant where
$G$ acts on itself by left multiplication. The map $\iota$ is
obviously $G$-invariant.

2) The local system $\cE$ has a canonical $G$-invariant structure: if
$x \in \fg$ and $\varphi \in \fg^*$, then $\cE_{(x,\varphi)} =
\wt{\varphi}_x$ and thus for all $g \in G$, one has
$$
\cE_{g(x,\varphi)} = \cE_{(gxg^{-1},{}^g\varphi)} =
\wt{({}^g\varphi)}_{gxg^{-1}} = \wt{\varphi}_{x}
$$

\noindent It follows from (1) and (2) that all the sheaves appearing
in the proof has a $G$-equivariant structure. It is easy to see using
(2) that the $G$-equivariant structure on $V = \iota^*\wt{u}^*(\cE)$
is the one that comes from the $G$-invariance of $v$ (i.e., $v(x,y) =
v(gxg^{-1},gy)$ for $g \in G$). Fix $g \in G$ and $x \in \fg$. Then
$g$ maps $V|_{x \times G}$ to $V|_{gxg^{-1} \times G}$ and thus one
gets an isomorphism $((\pr_1)_!V)_x\rightarrow
((\pr_1)_!V)_{gxg^{-1}}$. Let us put $g = x$. The resulting
automorphism of $((\pr_1)_!V)_x$ comes from the invariance of $v$
under $(x,y) \mapsto (x,xy)$ along the fiber $x \times G$ of
$\pr_1$. It follows that:

\begin {lem}
\label{automorphism}
The map $v$ is invariant under the automorphism $(x,y) \mapsto (x,xy)$
of $\fg \times G$ and the corresponding automorphism of $T$
{\em(}coming from the isomorphism of {\em Lemma \ref{inverseofl})} equals
$\theta_T$.
\end {lem}

\subsection {Proof  of Proposition $\ref{reformulation}$}
Since $\cF^{-1}$ is a quasi-inverse of $\cF$, one has the isomorphisms
$\ov{\ccL} \cong \cF(T)$ and thus $\ccL_\x \cong \cF(T)_\x$. 

\begin {lem}
\label {fd}
If $M \in \der(\fg)$, then $\cF(M)_\chi \cong R\Gamma_c(M \otimes
\wt{\x})$.
\end {lem}

\begin {proof}
Let $i$ be the inclusion of $\x$ in $\fg^*$ and let $j: \fg
\rightarrow \fg \times \fg^*$ be the inclusion $g \mapsto
(g,\x)$. One has the cartesian square:
$$
\xymatrix {
\fg \ar[r] \ar[d]^-{j} & \{\x\} \ar[d]^-{i}\\
\fg \times \fg^* \ar[r] & \fg^*
}
$$
Note that $j^*(\cE) = \wt{\x}$. Using the projection formula, one
obtains:
\begin {align*}
\cF(M)_\x &= (\pr'_!(\pr^*(M) \otimes \cE))_\x\\
&= i^*\pr'_!(\pr^*(M) \otimes \cE)\\
&= R\Gamma_c(j^*(\pr^*(M) \otimes \cE))\\
&= R\Gamma_c(M \otimes \wt{\x}),\ \text{qed} \hfill \qedhere
\end {align*}
\end {proof}
Let $i = n+d$. One has
\begin {align*}
T \otimes \wt{\x} &\cong ((\pr_1)_!(V) \otimes \wt{\x})[2i](i)\\
&= ((\pr_1)_!(V \otimes \pr_1^*(\wt{\x}))[2i](i)\\
&= ((\pr_1)_!(W))[2i](i)
\end {align*}
Thus it follows from the Lemma \ref{fd} that:
$$
\cF(T)_\x \cong R\Gamma_c((\pr_1)_!(W)[2i](i)) =
R\Gamma_c(W)[2i](i)
$$
As $\ccL_\x \cong \cF(T)_\x$, we get an isomorphism
$$
\sigma: \ccL_\x \rightarrow V = H^{2i}_c(W)
$$
(we can ignore the Tate twist since $k$ is algebraically
closed). Further, under the identification of $\ccL_\x$ and $\cF(T)_\x$,
the automorphism $\theta'_\ccL$ of $\ccL_\x$ is mapped to
the restriction of $\cF(\theta_T): \cF(T) \rightarrow \cF(T)$ to the
stalk $\cF(T)_\chi$. This automotphism is obtained by applying
$R\Gamma_c$ to $\theta_T \tensor \id_{\wt{\x}}: T \otimes \wt{\x}
\rightarrow T \otimes \wt{\x}$. Using the description of $\theta_T$
from Lemma \ref{automorphism}, it follows that $\theta_T \tensor
\id_{\wt{\x}}$ is the one induced by the invariance of $W$ under the
automorphism $(x,y) \mapsto (x,xy)$ of $\fg \times G$. Thus $\sigma$
maps $\theta'_\ccL$ to $\eta$. Finally, it remains to prove that
$\sigma$ is 
$\Gamma$ invariant. This follows easily from the
$G_\x$-equivariant structure of $V$ described before Lemma
\ref{automorphism}. The details are omitted.

\section {Quasi-polarizations}
In this section we introduce the notion of quasi-polarizations. This
will be used in computing the $\Gamma$ module structure on $V$ as well
as the automorhphism $\eta \in \Aut(V)$.
\subsection {Isotropic subgroups}
In this section the we shall only use the commutative group structure
on $\fg$. The results below apply equally well whenever $E$ is any
skew-symmetric biextension of $\fg$. 

In the following, we fix an additive {\em connected} subgroup $\fh$ of
$\fg$. Let $\fh^\perp$ denote the orthogonal complement of $\fh$ with
respect to $\E$. The group $\fh^\perp$ consists of all $x \in \fg$
such that the restriction of $\E$ to $x \times \fh$ 
is trivial. If $f: \fg \rightarrow \fg^*$ is the homomorphism
corresponding to $\E$, then $\fh^\perp$ is the kernel of the
composition: 
$$
\fg \xrightarrow{f} \fg^* \project \fh^*
$$
where the second map is the dual of the inclusion $\fh \hookrightarrow
\fg$. Note that $\ker \E = \fg^\perp = \ker f$. For notational
clarity, we shall denote $\ker E$ by $\fc$ (recall that $\ker E =
\fg_\x$ by Lemma \ref{kernel}). In the following, we assume that
$\fc^\circ \subset \fh$ and $\fh$ is {\em isotropic}, i.e., $\fh \subset
\fh^{\perp}$.

\begin {lem}
\label {dimension}
$\dim \fh^\perp + \dim \fh = \dim \fg + \dim \ker \E$
\end {lem}
\begin {proof}
We claim that the image of $f$ is the subgroup $(\fg/\fc^\circ)^*
\subset \fg$ consisiting of all elements of $\fg^*$ whose restriction
to $\fc^\circ$ is trivial. Indeed, it is clear that $\im f \subset
(\fg/\fc^\circ)^*$. Further, the kernel of $f$ is $\fc$, whence  there
exists an injection $f^\p:\fg/\fc \inject (\fg/\fc^\circ)^*$ such that
$f$ factors through $f^\p$, i.e., $f$ is the composition: 
$$
\fg \project \fg/\fc \xarrow{\inject}{f^\p} (\fg/\fc^\circ)^* \inject \fg^*
$$
The map $f^\p$ is an injective group homomorphism and thus a closed
immersion. Further, one has: 
$$
\dim (\fg/\fc^\circ)^* = \dim \fg/\fc^\circ = \dim \fg/\fc
$$
whence $f^\p$ is an isomorphism. Thus $\im f = (\fg/\fc^\circ)^*$, as
claimed. Let $\varphi$ denote the composition $\fg \xrightarrow{f}
\fg^* \surject \fh^*$. By definition, $\fh^\perp$ is the kernel of
$\varphi$. Further, the image of $\varphi$ is the image of
$(\fg/\fc^\circ)^*$ under the projection $\fg^* \project \fh^*$. As
$(\fg/\fc^\circ)^*$ surjects onto $(\fh/\fc^\circ)^*$ (since its dual
is injective), it follows that $\im \varphi =
(\fh/\fc^\circ)^*$. Therefore $\fg/\fh^\perp$ is isomorphic to
$(\fh/\fc^\circ)^*$. The equality now follows as $\dim
(\fh/\fc^\circ)^* = \dim \fh/\fc^\circ = \dim \fh - \dim \fc^\circ$. 
\end {proof}

Let $\fh^\p$ denote the neutral connected component of $\fh^\perp$.

\begin {lem}
\label {surjective}
The map $\pi_0(\fc) \rightarrow \pi_0(\fh^\perp)$
induced by the inclusion $\fc \subset \fh^\perp$ is surjective.
\end {lem}

\begin {proof}
Let $x \in \fh^\perp$. It suffices to find $y \in \fh^\p$ such that
$x-y \in \fc$, i.e., $f(x) = f(y)$. Note that $f(\fh^\perp)$, by
definition of $\fh^\perp$, is a subset of $(\fg/\fh)^*$. Further, we
claim that $f(\fh^\p) = (\fg/\fh)^*$. Indeed, if  $f^\p:\fh^\p
\rightarrow (\fg/\fh)^*$ is the map obtained by restricting $f$, then
it suffices to show that: 
$$
\dim \fh^\p = \dim (\fg/\fh)^* + \dim \ker f^\p
$$
But $\ker f^\p = \fh^\p \cap \fc$, which has the same dimension as
$\fc$ (because $\fh^\p$ contains $\fc^\circ$). Thus the above equality
follows from the Lemma. As $f(\fh^\p) = (\fg/\fh)^*$, there exists $y
\in \fh^\p$ such that $f(x) = f(y)$, qed. 
\end {proof}

As $\fh$ is connected, we have $\fh \subset \fh^\p$. One therefore has
the following lattice of subgroups of $\fg$: 

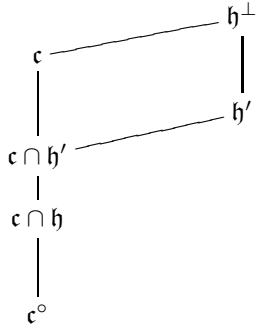
$$
\xymatrixrowsep{.1pc}
\xymatrix {       									&     & \fh^\perp  \ar@{-}[lld] \ar@{-}[ddd] &     & \\
\fc \ar@{-}[ddd]            &     &\\
&     &                     &\\
&     & \fh^\p \ar@{-}[lld] \\
\fc \cap \fh^\p \ar@{-}[dd] &     &\\
&&\\
\fc \cap \fh \ar@{-}[dddd]  &     &\\
&     &\\
&     &                     &     & \\ 
&     &                     &     & \\ 
\fc^\circ &    &\\
}
$$
The inclusions along the vertical segments have finite index.
\begin {cor}
\label {parallelogram}
The homomorphism $\fc/(\fc\, \cap\, \fh') \rightarrow \fh^\perp/\fh'$
induced by the inclusion $\fc \inject \fh^\perp$ is an isomorphism.
\end {cor}
\begin {proof}
Injectivity is obvious and surjectivity follows from Lemma
\ref{surjective}. 
\end {proof}

Recall that $\pi_0(\fc)$ was denoted by $\fp$. Define a filtration: 
$$
0 \subset \fp_\fh \subset \fp_{\fh^\p} \subset \fp
$$
of $\fp$ by setting $\fp_\fh = (\fc \cap \fh)/\fc^\circ$,
$\fp_{\fh^\p} = (\fc \cap \fh^\p)/\fc^\circ$. 

\begin {lem}
\label {orthogonal}
$\fp_\fh$ is isotropic with respect to the quadratic form $q$ on
$\fp$. Further, one has $\fp_\fh^\perp = \fp_{\fh^\p}$. 
\end {lem}

\begin {proof}
Since $\fh$ is isotropic, the restriction of $E$ to $\fh \times \fh$
gives the trivial biextension of $\fh \times \fh$. This implies that
$q|_{\fp_\fh} = 0$, whence the first assertion. Let $B: \fp \times \fp
\rightarrow \QZp$ be the biadditive map induced by $q$. Since the
restriction of $E$ to $\fh \times \fh^\p$ is trivial, it follows that
$B(\fp_\fh,\fp_{\fh^\p}) = 0$, i.e., $\fp_{\fh^\p} \subset
\fp_\fh^\perp$. It suffices to show that $\Card(\fp_{\fh'}) =
\Card(\fp/\fp_\fh)$, i.e., the cadinalities of $(\fc \cap
\fh)/\fc^\circ$ and $\fc/(\fc \cap \fh')$ are equal. Consider the
biextension of $\fh \times \fg$ obtained from $E$ by 
restriction. It induces a nondegenerate pairing:
$$
\frac{\fc \cap \fh}{\fc^\circ} \times \frac{\fh^\perp}{\fh'}
\rightarrow \QZp
$$
In particular, one has
$
\Card((\fc\cap \fh)/\fc^\circ) = \Card(\fh^\perp/\fh')
$.
The desired equality now follows from Corollary \ref{parallelogram}.
\end {proof}

\subsection {Quasi-polarizations: definition and existence results}

In this section, the hypothesis on the nilpotence class of $\fg$ is
not used. We retain the notations introduced in the previous section. 
\begin {df}
A {\em quasi-polarization} of $\E$ is a {\em connected, Lie
  subring} $\fh$ of $\fg$ such that 
\vspace {-6pt}
\begin {enumerate}[(i)]
\setlength {\itemsep} {-4pt}
\item $\fc^\circ \subset \fh$ and $\fh$ is isotropic;
\item $\fh^\perp$ is a Lie ring.
\end {enumerate}
\end {df}

A quasi-polarization $\fh$ is called a \textbf{Heisenberg
  polarization} if $[\fh^\perp, \fh^\p] \subset \fh$.

\begin {prop}
\label {heisenbergpolarization}
A Heisenberg polarization always exists.
\end {prop}

To prove this, one uses  the following:

\begin {prop}
\label {biggerpolarization}
Let $\fh$ be a quasi-polarization. If $\fh$ is not Heisenberg, then
there exists a quasi-polarization $\fr$ such that $\fr$ strictly
contains $\fh$ and
$
[\fh^\perp,\fh^\perp] \subset \fr^\perp \subset \fh^\perp
$.
\end {prop}

The proof of this is given in Section
4.3. To deduce Proposition \ref{heisenbergpolarization}, one starts with
the quasi-polarization $\fh = \fc^\circ$ and applies Proposition
\ref{biggerpolarization} to get an increasing chain of
quasi-polarizations terminating in a Heisenberg polarization. For
Heisenberg polarizations, one has the following analogue of
Proposition \ref{biggerpolarization}:

\begin {lem}
\label {biggerheisenbergpolarization} 
Let $\fh$ be a Heisenberg polarization and let $\fr$ be any
connected isotropic additive subgroup such that $\fh \subset \fr \subset
\fh^\p$. Then $\fr$ is a Heisenberg polarization and $
[\fh^\perp, \fh^\perp] \subset \fr^\perp \subset \fh^\perp
$.
\end {lem}

\begin {proof}
Indeed, since $[\fh^\p,\fh^\p] \subset \fh$, any additive subgroup of
$\fh^\p$ containing $\fh$ is a Lie ring. Thus $\fr$ is a Lie ring. It
remains to verify that $[\fh^\perp, \fh^\perp] \subset \fr^\perp
\subset \fh^\perp$ (this automatically implies that $\fr^\perp$ is a
Lie ring). The second inclusion is clear since $\fh \subset \fr$. For
the first inclusion, let us fix $x$, $y \in \fh^\perp$ and let $z =
[x,y]$. We need to show that the restriction $E|_{z \times \fr}$ is
trivial. Let $k \in \fr$. The Jacobi identity
$$
[z,k] = [x,[y,k]] + [y,[k,x]]
$$
implies that $E|_{z \times \fr}$ is the sum of two $\QZp$ torsors
$E^1, E^2$ where:
$$
E^1_{(z,k)} = \x_{[x,[y,k]]}, \quad E^2_{(z,k)} = \x_{[y,[k,x]]}
$$
It suffices to check that $E^1$, $E^2$ are trivial. Note that $[y,k]
\in \fh$ (since $[\fh^\perp,\fh'] \subset \fh$). Thus $E^1$ is
obtained by pulling back $E|_{x \times \fh}$ to $z \times \fr$ via the
map  $(z,k) \mapsto (x,[y,k])$. As $x \in \fh^\perp$, it follows that
$E_{x \times \fh}$ is trivial, whence so is $E^1$. One shows similarly
that $E^2$ is trivial.
\end {proof}

Let $\rank E$ denote the integer $\dim \fg - \dim \fc$.
\begin {cor}
\label{optimalpolarization}
If $\rank E$ is even, then there exists a quasi-polarization $\fr$
such that 
$\dim \fr = \dim \fr^\perp$ {\em(}equivalently, $\fr = \fr'${\em)}.
\end {cor}

One can obtain an additive version of this corollary (disregarding the
Lie ring structures) by replacing ``polarization $\fr$" by
``connected, isotropic (additive) subgroup $\fr$ containing $(\ker
\E)^\circ$". In this case, the result is proved in
Proposition A.28 of \cite{mitya}. 

\begin {proof}[Proof of Corollary $\ref{optimalpolarization}$]
Let $\fh$ be a Heisenberg polarization and let $\E^\p$ be the
biextension of $\fh^\p \times \fh^\p$ obtained from $\E$ by
restriction. We claim that $\rank E'$ is even. Indeed, the
homomorphism $f':\fh' \rightarrow \fh^{\p*}$ corresponding to $E'$ is 
the composition 
$$
\fh' \inject \fg \xrightarrow{f} \fg^* \rightarrow \fh^{\p*}
$$
Thus $\ker E' = \ker f' = \fh' \cap {\fh^\p}^\perp$. We claim that
$(\ker E')^\circ = \fh$. It suffices to show that $\fh$ is neutral
connected component of $\fh^{\p\perp}$. Since $\fh' \subset
\fh^\perp$, the restriction of $E$ to $\fh \times \fh'$ is
trivial. Thus $\fh \subset \fh^{\p\perp}$. Hence it remains to show
that $\dim \fh^{\p\perp} = \dim \fh$. Indeed, replacing $\fh$ by $\fh'$ in
Lemma \ref{dimension} gives $\dim \fh^{\p\perp} = \dim \fg + \dim \fc
- \dim \fh'$, which equals to $\dim \fh$ (again by Lemma
\ref{dimension}). Therefore $\rank E' = \dim \fh' - \dim \fh$, whence
it follows from Lemma \ref{dimension} that $\rank E - \rank E' =
2(\dim \fh - \dim \fc)$, qed. Applying Proposition A.28 of
\cite{mitya} to $E'$, one obtains a additive subgroup $\fr$ of $\fh'$
such that $\fh \subset \fr$ and $\dim \fr = \dim \fr^{\perp_\p}$,
where $\fr^{\perp_\p}$ is the orthogonal complement of $\fr$ in
$E^\p$. But $\fr^{\perp_\p} = \fr^\perp \cap \fh'$, whence $\dim \fr =
\dim \fr^\perp$. It follows from
Lemma \ref{biggerheisenbergpolarization} that $\fr$ is a
polarization. This completes the proof.
\end {proof}

\subsection {Proof of Proposition $\ref{biggerpolarization}$}
We shall need the elementary result:
\begin {lem}
\label {commutator}
Let $\fa, \fb$ be Lie subrings of $\fg$ such that $\fa$ is
connected. Then $[\fa,\fb]$ is also connected. 
\end {lem}
\begin {proof}
Let $\ft$ be any Lie ring. The image of the commutator map $\ft^\circ
\times \ft^\circ \rightarrow \ft$ is connected and contains 0, and
thus is contained in $\ft^\circ$. This shows that $\ft^\circ$ is also
a Lie ring. Now let $\ft = [\fa, \fb]$ and let $\fb_1$ be an arbitrary
connected component of $\fb$. The image of the commutator map $\fa
\times \fb_1 \rightarrow \ft$, for reasons similar to above, is
contained in $\ft^\circ$.  Therefore  whenever $a \in \fa, b \in \fb$,
one has $[a,b] \in \ft^\circ$. Since $\ft^\circ$ is a Lie ring, this
proves $\ft = \ft^\circ$, qed. 
\end {proof}

It follows that $\fh^\p$ is an ideal in $\fh^\perp$. Define a sequence
of Lie subrings $\fh^{(n)}$ of $\fh^\perp$ for integers $n \geq 0$ by:
\begin {align*}
  \fh^{(0)} &= \fh^\p  \\
  \fh^{(n)} &= [\fh^\perp, \fh^{(n-1)}] \qquad \text {for $n \geq 1$}
\end {align*}
One easily checks using Lemma \ref{commutator}  that the $\fh^{(n)}$
form a decreasing sequence of connected ideals of
$\fh^\perp$. Let $k$ be the integer such that $\fh^{(k)} \not \subset
\fh$ and $\fh^{(k+1)} 
\subset \fh$. Note that $k \geq 1$ since $\fh$ is not
Heisenberg. Define $\fr$ by 
$$
\fr = \fh + \fh^{(k)}
$$
It suffices to show that $\fr$ has the required properties. One
proceeds in various steps.

i) {\em $\fr$ is a connected Lie subring of $\fh^\perp$.} Indeed, both
$\fh$ and $\fh^{(k)}$ are connected, whence so is $\fr$. Further,
$\fh^{(k)}$ is an ideal in $\fh^\perp$ and thus $\fr$ is a Lie subring of
$\fh^\perp$.

ii) {\em If $\fa, \fb$ are two connected subgroups of $\fg$, then
  $(\fa + \fb)^\perp$ contains the neutral connected component of $\fa^\perp
  \cap \fb^\perp$.}  Since $(\fa + \fb)^\perp$ is contained in
$\fa^\perp \cap \fb^\perp$, it suffices to show that it has finite
index in $\fa^\perp \cap \fb^\perp$. Note that $(\fa + \fb)^\perp$ is
the kernel of the map $u: \fg \rightarrow (\fa + \fb)^*$ obtained as
the composition $$\fg \xrightarrow{f} \fg^* \rightarrow
(\fa+\fb)^*$$ where the second map sends $\varphi \in \fg^*$ to
$\varphi|_{\fa + \fb}$. Further, let $\pi:(\fa + \fb)^* \rightarrow
\fa^* \oplus \fb^*$ be the map $\varphi \mapsto
(\varphi_\fa,\varphi_\fb)$. Then $\fa^\perp \cap \fb^\perp$ is the
kernel of $\pi \circ u:\fg \rightarrow \fa^* \oplus \fb^*$. Thus it
suffices to show that $\pi$ has finite kernel, or equivalently,
$\pi^*$ is surjective (the equivalence easily follows from properties
(iii), (iv) of Section 2.3). This follows since $\pi^*: \fa \oplus \fb
\rightarrow \fa + \fb$ is the addition map.

iii) {\em $\fh^{(1)} = [\fh^\perp,\fh']$ is contained in $\fr^\perp$.} 
 As $\fh^{(1)}$ is connected, it suffices to show (using (ii)) that
 $\fh^{(1)} \subset \fh^\perp \cap 
 {\fh^{(k)}}^\perp$, i.e., $\fh^{(1)} \subset {\fh^{(k)}}^\perp$. This
 can be shown using the Jacobi identity in the style of the proof of
 Lemma \ref{biggerheisenbergpolarization}. One needs to use the fact
 that $[\fh^\perp, \fh^{(k)}] \subset \fh$.

iv) {\em $[\fh^\perp,\fh^\perp] \subset \fr^\perp \subset \fh^\perp$.}
The second inclusion is obvious since $\fh \subset \fr$. For the first
inclusion, note that $\fh^\perp = \fh' + \fc$ by Lemma
\ref{surjective}. As $\fc \in \fr^\perp$, it suffices to show that
$[\fc,\fh']$ and $[\fh',\fh']$ are contained in $\fr^\perp$. This
follows from (iii).

v) {\em $\fr^\perp$ is a Lie ring and $\fr$ is isotropic.} The first
statement is immediate from (iv). To show $\fr$ is isotropic, it
suffices to prove that $\fh + \fh^{(1)} \subset \fr^\perp$. Using
(iii), it remains to show that $\fh \subset \fr^\perp$. Indeed, since
$\fr \subset  \fh^\perp$,
the restriction of $E$ to $\fh \times \fr$ is trivial and thus $\fh
\subset {\fr}^\perp$.

This completes the proof of the proposition.

\section {Proof of Theorem $\ref{ribbon}$ when $\dim(\Omega)$ is even}
In this section we prove Theorem \ref{ribbon} when $\dim(\Omega)$ is even
by direct computation. Note that $$\dim(\Omega) = \dim(G) - \dim(G_\x)
= \dim(\fg) - \dim(\fc) = \rank(E)$$ Thus we can apply Corollary
\ref{optimalpolarization} in this case. We shall use the notation of Section 3.3 as
well as the following:\\ 

1) $\wt{E}$ denotes the local system on $\fg \times
\fg$ obtained by pulling back $\wt{\chi}$ using the commutator map
$\fg \times \fg \rightarrow \fg$ which sends $(u,v)$ to $[u,v]$. Note
that $\wt{E}$ is obtained from $E$ by means of $\psi:\QZp\inject
\ql^\times$. 

2) Let $X$ be a scheme over $k$ and let $F \in \der(X)$. We denote
$R\Gamma_c(X,F)$ by $\int_X F$.

3) If $Y$ is a locally closed subscheme of $X$ and $i: Y \inject X$ is
the inclusion, then we denote $i^*(F)$ by $F|_Y$ and
$R\Gamma_c(Y,F|_Y)$ by $\int_Y F$. \\

We note that if $Y$ is closed in $X$, then there is a natural morphism
$\int_X F 
\rightarrow \int_Y F$, which is an isomorphism if $\int_{X\bs Y} F$ is
zero (this is standard and follows from the distinguished triangle
described in Lemma 5.8 of \cite{mitya}).

\subsection {A basis of $V$}
For notational convenience, we shall identify $G$ with $\fg$ and
consider $W$ as a local system on $\fg \times \fg$.
Let $\fh$, $\fr$ be quasi-polarizations of $E$ such that $\fh \subset
\fr$ (and thus $\fr^\perp \subset \fh^\perp$).

\begin {prop}
\label {reduction}
If $[\fh^\perp,\fh^\perp] \subset \fr^\perp$, then the natural map
$$
\int_{\fh^\perp \times \fh^\perp} W \rightarrow \int_{\fr^\perp \times
  \fr^\perp} W
$$
is an isomorphism.
\end {prop}
The proof of this proposition is given in Section 5.4. By the results
of Section 4.2, it follows that there exists a sequence
$$
\fc^\circ = \fh_0 \subset \cdots \subset \fh_n = \fh
$$
of quasi-polarizations such that $\dim \fh = \dim \fh^\perp$ (i.e.,
$\fh = \fh'$) and $[\fh^\perp_i, \fh^\perp_i] \subset \fh^\perp_{i+1}$
for all $0 \leq i < n$. Thus Proposition \ref{reduction} implies that
the natural map $\int_{\fg \times \fg} W \rightarrow
\int_{\fh^\perp\times \fh^\perp}W$ is an isomorphism.
\begin {lem}
\label {trivial}
The local system $W|_{\fh^\perp\times \fh^\perp}$ is trivial.
\end {lem}
Before proving the lemma, we introduce some notation. One has:
\begin {align*}
w(x,y) = x -y^{-1}xy &= [y,x] - \frac{1}{2}[y,[y,x]] + \cdots\\
& = [y,\Phi(x,y)]
\end {align*}
where
$$
\Phi(x,y) = x - \frac{1}{2}[y,x] + \cdots +
\frac{(-1)^n}{(n+1)!}\ad^n(y)(x) + \cdots
$$
Let $\ld$ be the
automorphism of $\fg \times \fg$ which sends $(x,y)$ to
$(y,\Phi(x,y))$. Note that $W = \ld^*(\wt{E})$.
\begin {proof}[Proof of Lemma $\ref{trivial}$]
As $\ld$ maps $\fh^\perp \times \fh^\perp$ to itself, it suffices to
show that $\wt{E}|_{\fh^\perp \times \fh^\perp}$ is
trivial. Since $\fh$ is  isotropic, the biextension $E|_{\fh \times
  \fh}$ of $\fh \times \fh$ is trivial and thus $\wt{E}|_{\fh
  \times \fh}$ is a trivial local system. Let $\alpha$, $\beta \in
\pi_\circ(\fh^\perp)$ and let 
$\fh_\alpha$, $\fh_\beta$ be the corresponding connected components of
$\fh^\perp$. Choose $\wt{\alpha} \in \fh_\alpha$, $\wt{\beta} \in
\fh_\beta$. Let $x$, $y \in \fh$. Then the stalk of $\wt{E}$ at $(x +
\wt{\alpha}, y + \wt{\beta})$ is
\begin {align*}
\wx_{[x + \wt{\alpha}, y + \wt{\beta}]} &= \wx_{[x,y]} \otimes
\wx_{[x,\wt{\beta}]} \otimes \wx_{[\wt{\alpha},y]} \otimes
\wx_{[\wt{\alpha}, \wt{\beta}]}\\
&= \wt{E}_{(x,y)}\otimes \wt{E}_{(x,\wt{\beta})} \otimes
\wt{E}_{(\wt{\alpha},y)} \otimes \wt{E}_{(\wt{\alpha},\wt{\beta})}
\end {align*}
where the first equality uses the fact that $\wx$ is a multiplicative
local system. It remains to note that both $\wt{\alpha}$,
$\wt{\beta}$ belong to $\fh^\perp$, whence $\wt{E}|_{(\fh,\wt{\beta})}$
and $\wt{E}|_{(\wt{\alpha},\fh)}$ are trivial local systems, qed.
\end {proof}
Recall that we denoted $\pi_\circ(\fc)$ by $\fp$ (cf. Section
3.2). Let $\fb = \pi_\circ(\fh^\perp) = \fh^\perp/\fh$. By Lemma
\ref{surjective}, the map $\fp \rightarrow \fb$ induced by the
inclusion $\fc \inject \fh^\perp$ is surjective. Let $\fa$ be the
kernel of this map, so that one has an exact sequence
$$
0 \rightarrow \fa \rightarrow \fp \rightarrow \fb \rightarrow 0
$$
Since $\fh = \fh'$, it follows from Lemma \ref{orthogonal} that $\fa =
\fa^\perp$. Thus $\fa$ is a {\em Lagrangian} ideal (it can be shown
that it is abelian). In particular, one has $\Card(\fa) =
\Card(\fb) = \Card(\fp)^{1/2}$.

Fix a pair $(\alpha,\beta) \in \fb \times \fb$. Let $\wt{\beta}$ be an
element of $\fh_\beta$ such that $\wt{\beta} \in \fc$ (this is
possible since $\fp$ surjects onto $\fb$). For each such $\beta$, one
can choose a trivialization of $\fh_\alpha \times \fh_\beta$, as
follows. First, let $y_0$ be an element of $\fc$. Note that $W|_{\fg
  \times y_0}$ is a trivial local system. This follows since $\ld$
sends $\fg \times y_0$ to $y_0 \times \fg$ and $\wt{E}|_{y_0 \times
  \fg}$ is trivial. Further, the stalk of $W$ at $(0,y_0)$ is
$\wx_0$, which can be canonically identified with $\ql$ (since $\wx$
is multiplicative, one has an isomorphism $\wx_0 \otimes \wx_0 \cong
\wx_0$, which gives the identification). Thus $W|_{\fg \times y_0}$
has a natural trivialization. This holds in particular for $y_0 =
\wt{\beta}$. This trivialization induces trivializations
$$
t^\alpha_{\wt{\beta}}: \ql|_{\fh_\alpha \times \fh_\beta} \cong W|_{\fh_\alpha\times\fh_\beta}
$$
for each $\alpha \in \fb$ (the difference between the various
trivializations $t^\alpha_{\wt{\beta}}$ when $\wt{\beta}$ varies in
$\fh_\beta \cap \fc$ is computed in Lemma \ref{comparison} below). Let us
fix a map $s:\fb  \rightarrow \fc$
such that $s(\beta) \in \fh_\beta \cap \fc$ (in other words, $s$ is a
section of the composite map $\fc \rightarrow \fp \rightarrow
\fb$). We further assume that $s(0) = 0$. Using $t^\alpha_{s(\beta)}$,
one can identify $\int_{\fh_\alpha \times \fh_\beta}W$ with
$\int_{\fh_\alpha \times \fh_\beta} \ql$, which can be {\em canonically}
identified with $\ql[-2i](-i)$, where $i = \dim \fh_\alpha + \dim
\fh_\beta$ (note that $i = \dim \fg + \dim \fc$ by Lemma
\ref{dimension}). Let $1_{\alpha,\beta}$ denote the image of the
element $1 \in \ql = \ql(-i)$ (we ignore the Tate twist) in
$H^{2i}_c(\fh_\alpha \times \fh_\beta, W)$ by
$t^\alpha_{s(\beta)}$ . Then $\{1_{\alpha,\beta}\}$ form a
$\ql$-basis of $V$.

\subsection {The action of $\Gamma$ on $V$ and the autormorphism $\eta$}
Let $g$ be a point in $\fc = \Log(G_\x)$ and let $\sigma_g$ be the automorphism of
$\fg \times \fg$ which sends $(x,y)$ to $(gxg^{-1},gy)$. Since $W$ has
a $G_\x$-equivariant structure, one has isomorphisms $W_{(x,y)}
\rightarrow W_{\sigma_g(x,y)}$ for all $x$, $y \in G$. By transport of
structure, $\sigma_g$ sends a trivialization $t$ of
$W|_{\fh_\alpha\times\fh_\beta}$ to a trivialization of $W|_{\fh_{g\alpha
    g^{-1}}\times \fh_{g\beta}}$, which we denote by $\sigma_g(t)$.

\begin {lem} 
\label {transformation}
Let $\wt{\beta} \in \fh_\beta \cap \fc$. Then
  $\sigma_g(t^\alpha_{\wt{\beta}}) = t^{g\alpha
    g^{-1}}_{g\wt{\beta}}$.
\end {lem}
\begin {proof}
Note that $\sigma_g$ sends $\fg\times \wt{\beta}$ to $\fg \times
g\wt{\beta}$. Both $W|_{\fg \times \wt{\beta}}$ and $W|_{\fg \times
  g\wt{\beta}}$ are trivial with natural trivializations, coming from
the equalities $W_{0,\wt{\beta}} = \wx_0 = \ql$ and
$W_{(0,g\wt{\beta})} = \wx_0 = \ql$. It remains to show that the map
$W_{(0,\wt{\beta})} \rightarrow W_{(0,g\wt{\beta})}$ induced by $\sigma_g$ is
identity. This is immediate from the $G_\x$-equivariant structure of
$W$, defined in Section 3.4.
\end {proof}
Let $B: \fp \times \fp \rightarrow \QZp$ be the biadditive pairing
induced by $q$ and let $\wt{B}: \fp \times \fp \rightarrow \ql^\times$
be defined by $\wt{B} = \psi\circ B$. It is sometimes
convenient to extend the domain of $B$ and $\wt{B}$ to $\fc \times
\fc$ by means of the projection $\fc \times \fc \rightarrow \fp \times
\fp$.

\begin {lem}
\label {comparison}
Let $\wt{\beta}$, $\wt{\beta}' \in \fh_\beta \cap \fc$. Then
$t^\alpha_{\wt{\beta}'} = \wt{B}(\wt{\beta}'-\wt{\beta}, \Phi(\alpha,\beta))t^\alpha_{\wt{\beta}}$.
\end {lem}
The constant $\wt{B}(\wt{\beta}'-\wt{\beta},
\Phi(\alpha,\beta))$ in this formula is to be understood in the
following sense: since $\fa$ is isotropic, the pairing $\wt{B}$
induces a map from $\fa \times \fb$ to $\ql^\times$, which we again
denote by $\wt{B}$. Note that
$\Phi(\alpha,\beta) \in \fb$ and  $\wt{\beta}' - \wt{\beta}$
belongs to $\fh \cap \fc$, whence its image in $\fp$ is contained in
$\fa$. Thus $\wt{B}(\wt{\beta}'-\wt{\beta},
\Phi(\alpha,\beta))$ makes sense.

\begin {proof}[Proof of Lemma $\ref{comparison}$]
Using the automorphism $\ld$ of $\fg \times \fg$, we shall shift the
computation from the local system $W$ to $\wt{E}$. Note that:

1) $\ld$ sends $\fh_\alpha
\times \fh_\beta$ to $\fh_\beta \times
\fh_{\Phi(\alpha,\beta)}$. 

2) For each $y_0 \in \fc$, the
trivial local system $\wt{E}|_{y_0 \times \fg}$ admits a natural
trivialization coming from the equalities $\wt{E}_{(y_0,0)} = \wx_0 =
\ql$. Further, the
trivialization of $W|_{\fg \times y_0}$ described in Section 5.1 is
obtained by pulling back that of $\wt{E}|_{y_0 \times \fg}$
by $\ld$. 

Thus it remains to show that the trivializations of
$\wt{E}|_{\fh_{\beta}\times \fh_{\Phi(\alpha,\beta)}}$ induced by
the trivializations of $\wt{E}|_{\wt{\beta}'\times \fg}$ and
$\wt{E}|_{\wt{\beta}\times \fg}$ differ by the desired
constant. Choose $y_0 \in \fh_{\Phi(\alpha,\beta)} \cap \fc$. Then these
trivializations give isomorphisms $\ql \cong
\wt{E}_{(\wt{\beta}',y_0)}$ and $\ql \cong
\wt{E}|_{(\wt{\beta},y_0)}$, which induce trivializations $t_1$ and
$t_2$ of $\wt{E}|_{\fg \times y_0}$. Further, $\wt{E}|_{\fg \times y_0}$
is equipped with a natural trivialization $t:\ql|_{\fg \times y_0}
\rightarrow \wt{E}|_{\fg \times y_0}$ (defined as in (2)). It follows
from the definition of $B$ (cf. Section 3.1) that $t_1 =
\wt{B}(\wt{\beta}',y_0)t$ and $t_2 = \wt{B}(\wt{\beta},y_0)t$. Hence
$t_1$ and $t_2$ differ by $\wt{B}(\wt{\beta}'-\wt{\beta},y)$. This
completes the proof.
\end {proof}
We can now describe the action of $\Gamma$ on the basis elements
$1_{\alpha,\beta}$ of $V$.
\begin {cor}
\label {gammaaction}
Let $g \in \Gamma$. Then $g(1_{\alpha,\beta}) =
\wt{B}(gs(\beta)-s(g\beta),\Phi(g\alpha
g^{-1},g\beta))1_{g\alpha g^{-1},g\beta}$. 
\end {cor}

\begin {proof}
We may assume that $g \in \fc = \Log(G_\x)$. By Lemma
\ref{transformation}, the map $\sigma_g: \fh_\alpha \times \fh_\beta
\rightarrow \fh_{g\alpha g^{-1}} \times \fh_{g\beta}$ transforms the
trivialization $t^\alpha_{s(\beta)}$ to $t^{g\alpha
  g^{-1}}_{gs(\beta)}$. It remains to note that by Lemma
\ref{comparison}, the trivializations $t^{g\alpha g^{-1}}_{gs(\beta)}$
and $t^{g\alpha g^{-1}}_{s(g\beta)}$ differ by
$\wt{B}(gs(\beta)-s(g\beta), \Phi(g\alpha g^{-1},g\beta))$, qed.
\end {proof}
It remains to compute the automorphism $\eta$. Let $x_0$, $y_0 \in
\fc$. Then $W_{(x_0,y_0)} = \wx_{x_0 - y_0^{-1}x_0y_0} =
W_{(x_0,x_0y_0)}$. Further, the trivializations of $W|_{\fg \times
  y_0}$ and $W|_{\fg \times x_0y_0}$ give isomorphisms $t_x: \ql
\rightarrow W_{(x,y_0)}$ and $t'_x: \ql \rightarrow W_{(x,x_0y_0)}$
for each $x \in \fg$.
\begin {lem}
One has $t_{x_0} = {\wt{q}(x_0)}^{-1}t'_{x_0}$.
\end {lem}

\begin {proof}
For notational clarity, we shall denote the stalk $\wx_x$ by
$\wx(x)$. One has
\begin {align*}
W_{(x,x_0y_0)} &= \wx(x - y_0^{-1}x_0^{-1}xx_0y_0)\\ &= \wx(x -
x_0xx_0^{-1})\otimes \wx(x_0xx_0^{-1} - y_0^{-1}x_0^{-1}xx_0y_0)\\
&= W_{(x,x_0)}\otimes W_{(x_0^{-1}xx_0,y_0)}
\end {align*}
Note that $W|_{\fg \times x_0}$ is also trivial and let $t''_{x}:\ql
\rightarrow W_{(x,x_0)}$ be the trivialization. It follows from the
above that
$$
t'_x = t''_x \otimes t_{x_0xx_0^{-1}}
$$
Putting $x = x_0$, we get $t'_{x_0} = t''_{x_0}\otimes t'_{x_0}$. Thus
it remains to show that $t''_{x_0}:\ql \rightarrow W_{(x_0,x_0)} =
\wx(0) = \ql$ is multiplication by $\wt{q}(x_0)$. Note that the
automorphism $\ld$ of $\fg \times \fg$ takes $\fg \times x_0$ to
$x_0 \times \fg$ and the trivialization $t''$ to the (natural)
trivialization of $\wt{E}|_{x_0 \times \fg}$. Thus it remains to check
that the map $\ql \rightarrow \wt{E}_{(x_0,x_0)} = \ql$ coming from
the trivialization of $\wt{E}|_{x_0 \times \fg}$ is multiplication
by $\wt{q}(x_0)$. But this follows from the definition of the
quadratic form $q$ (cf. Section 3.1).
\end {proof}
Let $\tau$ be the automorphism of $\fg \times \fg$ which sends $(x,y)$
to $(x,xy)$. Let $\wt{\alpha} \in \fh_\alpha \cap \fc$ and $\wt{\beta}
\in \fh_\beta \cap \fc$. Putting $x_0 = \wt{\alpha}$ and $y_0 = \wt{\beta}$ in
the lemma shows that
\begin {cor}
\label{preetaaction}
$\tau$ takes the trivialization $t^{\wt{\beta}}_\alpha$ to the
trivialization
$\wt{q}(\wt{\alpha})^{-1}t^{\wt{\alpha}\wt{\beta}}_\alpha$ of
$\fh_\alpha \times \fh_{\alpha \beta}$.
\end {cor}
Using the corollary, it is easy to give a formula for
$\eta(1_{\alpha,\beta})$. However, we shall only need the following:
\begin {cor}
\label {etaaction}
$\eta(1_{\alpha,0}) = \wt{q}(s(\alpha))^{-1}1_{\alpha,\alpha}$.
\end {cor}
\begin {proof}
Putting $\wt{\alpha} = s(\alpha)$ and $\wt{\beta} = 0$ in
Corollary \ref{preetaaction}.
\end {proof} 
\subsection {Proof of Theorem \ref{ribbon}}
We shall prove the theorem by means of explicit computation using
Corollary \ref{gammaaction} and Corollary \ref{etaaction}. Put
$$
u = \sum_{\alpha \in \fb} 1_{\alpha,0}
$$
\begin {lem}
The elements $gu$, $g \in \Gamma$ form a $\ql$-basis of $V$.
\end {lem}

\begin {proof}
Let $\wt{\beta}$ denote the image of $s(\beta) \in \fc$
in $\fp$. Fix $\beta_0 \in \fb$. Let $\pi: \fp \rightarrow \fb$ denote
the projection. Note that $\pi^{-1}(\beta_0) = \fa +
\wt{\beta}_0$. Let $g \in \pi^{-1}(\beta_0)$. It follows from
Corollary \ref{gammaaction} that 
$$
g(1_{\alpha,0}) = \wt{B}(g - s(g), \Phi(g\alpha g^{-1},g))1_{g\alpha
  g^{-1}, g}
$$
The pair $(g\alpha g^{-1},g) \in \fb \times \fb$ is none other than $(\beta_0\alpha
\beta_0^{-1}, \beta_0)$. Thus $$\wt{B}((g - s(g), \Phi(g\alpha g^{-1},g))
= \wt{B}(g - \wt{\beta}_0, \Phi(\beta_0\alpha
\beta_0^{-1},\beta_0))$$ Writing $g = a + \wt{\beta}_0$, it follows
that $g(1_{\alpha,0}) = \wt{B}(a,\Phi(\beta_0\alpha\beta_0^{-1},\beta_0))
1_{\beta_0\alpha \beta_0^{-1}, \beta_0}$ and thus
\begin {align*}
g(u) &= \sum_{\alpha \in \fb}  \wt{B}(a,\Phi(\beta_0\alpha\beta_0^{-1},\beta_0))
1_{\beta_0\alpha \beta_0^{-1}, \beta_0}\\
&= \sum_{\alpha \in \fb}  \wt{B}(a,\Phi(\alpha,\beta_0))
1_{\alpha, \beta_0}
\end {align*}
The map $\alpha \rightarrow \wt{B}(a,\Phi(\alpha,\beta_0))$ from $\fb$ to
$\ql^\times$ is a character of $\fb$ (this follows since $\Phi(x,y)$ is
additive in $x$), which we denote by $\zeta_a$. Note that $\zeta_a$ is
trivial if and only if $a = 0$ (i.e., $g = \wt{\beta}_0$). Indeed, as
$\alpha$ varies over $\fb$, the elements $\Phi(\alpha,\beta_0)$
exhaust $\fb$. Thus if $\zeta_a = 0$, then $\wt{B}(a,\fb) = 0$, which
implies $a = 0$. Also, $\zeta_{a+b} = \zeta_a\zeta_b$, whence the map
$a \mapsto \zeta_a$  from $\fa$ to $\fb^*$ is an isomorphism of
groups. It follows 
that the elements $g(u)$ as $g$ varies over $\pi^{-1}(\beta_0)$ form a
basis of the subspace of $V$ generated by $\{1_{\alpha,\beta_0}\}$,
$\alpha \in \fb$. This completes the proof.
\end {proof}
In view of the lemma, it suffices to show that $\eta(u) =
\wh{q}(u)$. Let $h_{\beta_0}$ be the element of $\ql[\Gamma]$ which
sends $u$ to $1_{\beta_0,\beta_0}$. Now one has
$$
(a+\wt{\beta}_0)(u) = \sum_{\alpha \in \fb}
\zeta_a(\alpha)1_{\alpha,\beta_0}
$$
It follows from the Fourier inversion formula for finite abelian group
that
$$
h_{\beta_0} = \frac{1}{\Card(\fa)}\sum_{a \in
  \fa}\zeta_a(-\beta_0)(a + \wt{\beta}_0)
$$
Note that $\zeta_a(-\beta_0) = \wt{B}(a,\Phi(-\beta_0,\beta_0)) =
\wt{B}(a,-\beta_0)$. Further:
$$
\wt{B}(a,-\beta_0) = \frac{\wt{q}(a)\wt{q}(\wt{\beta}_0)}{\wt{q}(a +
  \wt{\beta}_0)} = \frac{\wt{q}(\wt{\beta}_0)}{\wt{q}(a +
  \wt{\beta}_0)}
$$
where the second equality follows from the fact that $\wt{q}(\fa) = 1$
(indeed, the biextension $E|_{\fh \times \fh}$ is trivial, whence the
result). Thus
$$
h_{\beta_0} = \frac{\wt{q}(\wt{\beta}_0)}{\Card(\fa)}\sum_{g \in
  \pi^{-1}(\beta_0)}\frac{1}{\wt{q}(g)}g
$$
As $\eta(u) = \sum_{\beta_0 \in \fb}
\wt{q}(\wt{\beta}_0)^{-1}1_{\beta_0,\beta_0}$, it suffices to show that
$\wh{q}$ is equal to
\begin {align*}
\sum_{\beta_0 \in \fb} \frac{1}{\wt{q}(\wt{\beta}_0)} h_{\beta_0}
&= \frac{1}{\Card(\fa)} \sum_{g \in \Gamma} \frac{1}{\wt{q}(g)}g
\end {align*}
Comparing with Lemma \ref{formulaforqhat}, it remains to show that
$G(\fp,\wt{q}) = \Card(\fa)$. This
is a consequence of Corollary 6.2, \cite{dgno}.

\subsection {Proof of Proposition \ref{reduction}}
Let $A = \fr^\perp \times \fr^\perp$, $A' = \fh^\perp \times
\fh^\perp$. It suffices to show that $\int_{A'\bs A} W = 0$. Let $B =
\fh^\perp\times \fr^\perp$. One has $A \subset B \subset A'$. Let 
\begin {align*}
U &= A'\bs B = \fh^\perp \times (\fh^\perp\bs \fr^\perp)\\
V &= B\bs A = (\fh^\perp\bs \fr^\perp) \times \fr^\perp
\end {align*}
It remains to prove that: i) $\int_U W = 0$ and ii) $\int_V W = 0$.




\begin {proof}[Proof of {\em (i)}.]
Consider the projection $\pr_2: U \rightarrow \fh^\perp\bs
\fr^\perp$. It suffices to show that $(\pr_2)!(W|_U) = 0$. If $y_0 \in
\fh^\perp\bs \fr^\perp$, then the stalk of $(\pr_2)_!(W|_U)$ at $y_0$
is $\int_{\fh^\perp \times y_0} W$. We shall show that this is
zero. Let $i: \fg \rightarrow \fg$ be the map $x \mapsto x -
y_0^{-1}xy_0$. Then $\int_{\fh^\perp \times y_0} W \cong
\int_{\fh^\perp} i^*(\wt{\x})$. As $i$ is additive, the local system
$i^*(\wt{\chi})$ is multiplicative. More precisely, one has
$$
i^*(\wx) = \wx \otimes ({}^{y_0}\wx)^{-1}
$$
We shall now use the following result:
\begin {lem}
\label {vanishing}
Let $\ccL$ be a multiplicative local system on a group scheme $H$. If
$\ccL|_{H^\circ}$ is nontrivial, then $R\Gamma_c(H,\ccL) = 0$.
\end {lem}
This is well known and is proved in Lemma 9.4, \cite{mitya}. Using
this, it remains to show that the local system $i^*(\wx)|_{\fh'}$ is
nontrivial, i.e., $\chi|_{\fh'}$ is not isomorphic to
${}^{y_0}\chi|_{\fh'}$. 
\begin {lem}
Let $\fa \subset \fg$ be a Lie subalgebra and let $b \in \fg$ be such
that $[b,\fa] \subset \fa$. Then $\chi|_\fa = {}^b\chi|_\fa$ if and
only if $b \in \fa^\perp$.
\end {lem}
\begin {proof}
This is a generalization of Lemma \ref{kernel} and the same proof
applies (with $x = b$ and $y \in\fa$).
\end {proof}
As $y_0 \in \fh^\perp$, one has $[y_0,\fh'] \subset \fh'$. Thus it
suffices to show that $y_0 \notin (\fh')^\perp$. This follows since
$y_0 \notin \fr^\perp$ and $(\fh')^\perp \subset \fr^\perp$ (as $\fr \subset
\fh'$). 
\end {proof}
\begin {proof} [Proof of {\em(ii)}]
Using the first projection $V \rightarrow \fh^\perp\bs
\fr^\perp$, it suffices to show that $\int_{x_0 \times \fr^\perp} W =
0$ for all $x_0 \in \fh^\perp\bs \fr^\perp$. Let $i: \fr^\perp \rightarrow
\fg$ be the map $y  \mapsto x_0 -
y^{-1}x_0y$. Then $\int_{x_0 \times \fr^\perp} W\cong \int_{\fr^\perp}
i^*(\wx)$. Unfortunately, $i$ is not additive. However, it suffices to
show that $\int_{y_0\fr}i^*(\wx) = 0$ for all $y_0 \in \fr^\perp$ (if
$\pi: \fr^\perp  \rightarrow \fr^\perp/\fr$ is
the projection map, then this would prove that
$\pi_!(i^*(\wx)) = 0$). Define $j:  \fr
\rightarrow \fg$ by $y \mapsto x_0 - y^{-1}y_0^{-1}x_0y_0y$. Then
$\int_{y_0\fr} i^*(\wx) \cong \int_\fr j^*(\wx)$. Put $t =
y_0^{-1}x_0y_0$. One has
$$
j(y) = x_0 - y^{-1}ty = (x_0 - t) + [y,t] - \frac{1}{2}[y,[y,t]] + \cdots
$$
Let $j_n$ denote that $n$-th term of this Lie series (i.e., $j_0 =
(x_0 - t)$, $j_1 = [y,t]$ etc.). Since $\wx$ is multiplicative, one
has
$$
j^*(\wx) = j_0^*(\wx) \otimes j_1^*(\wx) \otimes \cdots
$$
\begin {claim}
$j_n^*(\wx) \cong (\ql)_\fr$ if $n \neq 1$.
\end {claim}
\begin {proof}
Since $j_0$ is a constant map, this is obvious for $n = 0$. Suppose
that $n = 2$. Let $s: \fr \rightarrow \fg \times \fg$ be the map  $y
\mapsto (y,-\frac{1}{2}[y,t])$. Then $j_2^*(\wx) = s^*(\wt{E})$.
As $y \in \fr \subset \fh^\perp$ and $t \in \fh^\perp$, one has $[y,t] \subset
[\fh^\perp,\fh^\perp]$, which is contained in $\fr^\perp$ by
assumption. Thus $\im s \subset \fr \times \fr^\perp$. But
$\wt{E}|_{\fr \times \fr^\perp} = (\ql)_{\fr \times \fr^\perp}$
(cf. proof of Lemma \ref{trivial}),
whence $j_2^*(\wx) \cong (\ql)_\fr$. The proof for $n > 2$ is similar.
\end {proof}
It follows that $j^*(\wx) \cong j_1^*(\wx)$. Since $j_1$ is additive,
the local system $j_1^*(\wx)$ is multiplicative. If this is trivial,
then $t \in \fr^\perp$. As $t = y_0x_0y_0^{-1}$, this implies $x_0 \in
\fr^\perp$, which is absurd. It now follows from Lemma \ref{vanishing}
that $\int_\fr j^*(\wx) = 0$. This completes the proof (ii) as well as
the proof of the proposition.

\end {proof}

\end {document}